Mailbox
                 Nuovo:
                           Messaggio
                                        Evento
                                                 Attività
                                                          Nota

                                                                  Filtri
                                                                              Configurazioni
                                                                                                    Sito WEB

\newcount\notenumber

\def\note{\advance\notenumber by 1
\footnote{$^{(\the\notenumber)}$}}

\font\tenmsb=msbm10 \font\sevenmsb=msbm7 \font\fivemsb=msbm5
\newfam\msbfam
\textfont\msbfam=\tenmsb \scriptfont\msbfam=\sevenmsb
\scriptscriptfont\msbfam=\fivemsb \edef\msbfam{\ifcase\msbfam 0\or
1\or 2\or 3\or 4\or 5\or 6\or 7\or 8\or 9\or A\or B\or C\or D\or
E\or F\fi}

\ifx\sc\undefined
    \font\sc=cmcsc10   
\fi

\mathchardef\nmid="3\msbfam2D
\def\Bbb#1{{\fam\msbfam\relax#1}}
\def\N{{\Bbb N}}
\def\Q{{\Bbb Q}}

\def\Z{{\Bbb Z}}
\def\G{{\Bbb G}}

\def\C{{\Bbb C}}

\def\V{{\cal V}}

\def\1{{\bf C}}
\def\A{{\Bbb A}}

\def\h{{\bf h}}
\def\v{{\varphi }}
\def\u{{\bf u}}
\def\d{{\rm d}}

\def\w{{\bf w}}

\def\CVD{{\hfill\hfil{\lower 2pt\hbox{\vrule\vbox to 7pt
{\hrule width 6pt\vfill\hrule}\vrule}}}\par}


\font\title=cmr10 scaled 1200

\baselineskip=15pt
 
\hsize = 15truecm \vsize = 22truecm
 
\hoffset = 0.4truecm \voffset = 0.7truecm

 
\centerline{\title On composite lacunary polynomials and the proof
of a conjecture of Schinzel}\bigskip
 
\centerline{Umberto Zannier}\bigskip\bigskip
 
\noindent{\bf Abstract}.  Let $g(x)$ be a fixed non-constant
complex polynomial. It was conjectured by Schinzel that if
$g(h(x))$ has boundedly many terms, then $h(x)\in \C[x]$ must also
have boundedly many terms. Solving  an older conjecture raised by
R\'enyi and  by Erd\"os, Schinzel had proved this in the special
cases $g(x)=x^d$; however that method does not extend to the
general case. Here we prove the full Schinzel's conjecture
(actually in sharper form) by a completely different method.
Simultaneously we  establish an ``algorithmic" parametric
description of the general decomposition $f(x)=g(h(x))$, where $f$
is a polynomial with a given number of terms and $g,h$ are
arbitrary polynomials. As a corollary, this implies for instance
that a polynomial with $l$ terms and  given coefficients is
non-trivially decomposable if and only if the degree-vector lies
in the union of certain finitely many subgroups of $\Z^l$.
\bigskip
 
\centerline{*****}

\bigskip
 
 
\noindent{\sc Introduction}. The behaviour of (complex)
polynomials under the operation of composition  has been studied
by several authors, starting with J.F. Ritt (see [S2] for an
account of the theory). Here we deal with this aspect when some of
the involved polynomials are {\it lacunary} (also called {\it
sparse}), i.e. the number of their terms is viewed as fixed, while
the corresponding degrees (and coefficients) may vary. So, we
write $f(x)=a_1x^{m_1}+\ldots +a_lx^{m_l}$ for a lacunary
polynomial with (at most) $l$ terms and we study its {\it
decomposability}, i.e. the equation $f(x)=g(h(x))$, with
$g,h\in\C[x]$ of degree $>1$; both decomposable and lacunary
polynomials have played a special role in several (algebraical and
arithmetical) investigations (see e.g. [S2]).
 
 A trivial case occurs when $h(x)=ax^n+b$; now, $f$
is of the shape $g\circ h$ if and only if $n$ divides all the
degrees of the terms which occur in $f(x)$. For non-trivial
decompositions, in a recent paper we established a bound (which
will be useful later) for the degree of $g(x)$ (see Thm. 1 of
[Z]):\medskip
 
\noindent{\bf Theorem A}. ([Z], Thm. 1) {\it Suppose that
$g,h\in\C[x]$ are non-constant, that $h(x)$ is not of the shape
$ax^n+b$ and that $g(h(x))$ has at most $l$ terms. Then $\deg g\le
2l(l+1)$.}\medskip
 
This somewhat controls the polynomial $g(x)$. To control $h(x)$
 leads to subtler problems already in basic cases, like
$g(x)=x^2$:  it was conjectured by R\'enyi and independently by
Erd\"os in 1949 [E] that a bound for the number of terms of
$h(x)^2$  implies a bound for the number of terms of $h(x)$. In
1987 Schinzel [S] found an ingenious proof of this conjecture,
actually for all powers $h(x)^d$ (and he gave explicit bounds). He
went on to conjecture that for a fixed non-constant $g\in \C[x]$
such that $g(h(x))$ has at most $l$ terms, the number of terms of
$h(x)$ is bounded by a function only of $l$. (This generalized
conjecture, as we shall see, has   significant  implications in
the whole context.) He also remarked that his method  for the
powers $h(x)^d$  was insufficient for a general proof.
 
In this paper  we fully prove this conjecture of Schinzel,
actually in   sharper form, namely without fixing the polynomial
$g(x)$. We have:
\medskip
 
\noindent{\bf Theorem 1}. {\it There exists a (computable)
function $B$ on $\N$ such that if $g,h\in\C[x]$ are non-constant
polynomials and if $g(h(x))$ has at most $l$ terms, then $h(x)$
has at most $B(l)$ terms.}\medskip
 
Our arguments follow a completely different path with respect to
Schinzel's proof of the special case (so in particular they
provide an alternative proof of the R\'enyi-Erd\"os conjecture).
They mainly rely on a kind of modified Puiseux expansions   and on
a lower bound for approximations  by sums of $S$-units in function
fields (see Prop. 1 below); this may be viewed as a case of
Schmidt Subspace Theorem in function fields.
 
The present proofs would easily yield an explicit, though very
large,  estimation  for $B(l)$, but for simplicity we do not
calculate it here. \note{Schinzel [S] produces explicit bounds for
the special cases $g(x)=x^d$; we believe that the present method
leads to weaker bounds in those cases. For bounds in the opposite
direction see [E], [S2].}
\medskip

Theorem 1  in full generality  represents (together with Theorem
A)  an indispensable tool
 {\it to obtain  the classification of polynomials
$f(x)$, with at most $l$ terms, which are ``decomposable"}, i.e.,
of the shape $g(h(x))$ with $g,h$ of degree $>1$.  Simultaneously
with Theorem 1 we establish a complete ``algorithmic" description
in finite terms. That is, for any fixed $l$ we give an effective
procedure to write down a finite number of parametrizations for
all the equations $f(x)=a_1x^{m_1}+\ldots +a_lx^{m_l}=g(h(x))$
where by ``parametrization" we roughly mean:
 
  \centerline{``{\it algebraic variety for the coefficients-vector $+$ integer
lattice  for the degrees-vector}".}
 
We can rephrase this by saying  that we can  obtain all the
equations in question from finitely many ``generic equations" just
by substitution.  More precisely we have:\medskip
 
\noindent{\bf Theorem 2}. {\it Let $l$ be a positive integer.
There exist an integer $p$, finitely many affine varieties
$\V_j/\Q$, $j=1,\ldots ,J$,  and  polynomials $F_j,H_j\in
\Q[\V_j][z_1^{\pm 1},\ldots ,z_p^{\pm 1}]$, $G_j\in\Q[\V_j][z]$,
such that:
 
(i) $F_j=G_j\circ H_j$.
 
(ii)  $F_j$ has at most $l$ terms as a Laurent polynomial in
$z_1,\ldots ,z_p$ and $\deg_zG_j\le 2l(l+1)$.
 
(iii) If $f,g,h\in \C[x]$ are such that $f=g\circ h$, $h(x)\neq
ax^n+b$  and $f$ has at most $l$ terms, then, for some $j$ there
exist a point $P\in\V_j(\C)$ and integers $u_1,\ldots ,u_p$ such
that $f(x)=F_j(P,x^{u_1},\ldots ,x^{u_p})$, $g(x)=G_j(P,x)$,
$h(x)=H_j(P,x^{u_1},\ldots ,x^{u_p})$.
 
  Finally, one may effectively find
$p,J$, equations for the $\V_j$ and expressions for the
$F_j,G_j,H_j$.}\medskip
 
  See also the equivalent Theorem 2$^*$ below for an
alternative formulation.
 
Theorem 2 follows  rather easily from Theorems A, 1. However the
proof of Theorem 1 in turn involves a description like in Theorem
2, so in fact the proofs will appear at the same time.
 
We also note that Theorem 2 immediately implies for instance  the
following:\medskip
 
\noindent{\bf Corollary}.\ {\it For  $a_1,\ldots ,a_l\in\C$, there
exists a finite union $M=M(a_1,\ldots ,a_l)$ of  subgroups of
$\Z^l$ such that $a_1x^{m_1}+\ldots +a_lx^{m_l}\in\C[x]$ is
nontrivially decomposable  if and only if $(m_1,\ldots ,m_l)\in
M$.}
\medskip
 
We can also add that if $\Q(a_1,\ldots ,a_l)$ is finitely
presented the finitely many relevant subgroups are computable.
Moreover, a similar corollary holds concerning the decomposability
of $a_1x^{m_1}+\ldots +a_lx^{m_l}$ for some $(a_1,\ldots ,a_l)$
running through a given algebraic variety.

Some of the arguments should extend to Laurent polynomials (as in
[Z]), to rational functions and also to  equations of the form
$g(c_1x^{m_1},\ldots ,c_lx^{m_l},h(x))=0$, for a fixed $g\in
\C[X_1,\ldots ,X_l, Y]$, where  $c_i\in \C$, $m_i\in \N$ and where
$h(x)$ is a polynomial. In turn, this is related to a Bertini-type
theorem, for the irreducibility of the intersection of a
subvariety of $\G_m^n$ with families of algebraic subgroups or
cosets. Also in view of the fact that this topic falls somewhat
far from the present one, we do not treat it here.
\bigskip

\noindent{\sc Proofs}. We shall need a version of the Voloch and
Brownawell \& Masser ``$S$-unit equation theorem for function
fields". Actually, rather than $S$-unit equations we shall meet
approximations by $S$-units, and for our purposes the following
variant shall be useful, modelled on [Z2, Thm. 1]:\medskip
 
\noindent{\bf Proposition 1}. {\it Let $K/\C$  be a function field
in one variable, of genus $g$, and let $\v_1,\ldots ,\v_n\in K$ be
linearly independent over $\C$. Let $S$ be a finite set of  places
of $K$ containing all the poles of $\v_1,\ldots ,\v_n$ and also
all the zeros of $\v_{1},\ldots ,\v_r$. Further, put
$\sigma=\sum_{i=1}^n\v_i$. Then $$ \sum_{v\in
S}(v(\sigma)-\min_{i=1}^n v(\v_i))\le {n\choose 2}(\#
S+2g-2)+\sum_{i=r+1}^n\deg(\v_i). $$}\medskip
 
\noindent{\it Proof}. Following  [BM] and [Z2], for a non-constant
$t\in K$ and $\v_1,\ldots ,\v_n\in K$, we consider the Wronskian
$W_t(\v_1,\ldots ,\v_n)$, i.e. the determinant of the $n\times n$
   matrix whose $j$-th row-entries are the $(j-1)$-th derivatives of the $\v_i$'s with
respect to $t$.
  Since the $\v_i$ are
linearly independent over  $\C$, we have $W_t\neq 0$ by a
well-known criterion. Let $z\in K$ be another non-constant
element. Then we have the known, easily proved, formula
$
W_z(\v_1,\ldots ,\v_n)=({dt\over dz})^{n\choose 2}W_t(\v_1,\ldots
,\v_n).
$
For a place $v$ of $K$ we choose once and for all a   local
parameter $t_v$ at $v$ and we define $W_v:=W_{t_v}$. This depends
on the choice of $t_v$, but the   formula  shows that the order
$v(W_v)$ depends only on $v$.
 
Since $W_t\in K^*$, the formula also shows that
$\sum_vv(W_v)={n\choose 2}\sum_v(\d t/\d t_v)={n\choose 2}(2g-2)$.
 
  For $v\not\in S$ we have $v(\v_i)\ge 0$ for $i=1,\ldots ,n$, so $v(\d^m\v_i/\d
t_v^m)\ge 0$ for all $i,m$ and $v(W_v)\ge 0$.
 
For $v\in S$, let $j=j_v$ be an index such that
$v(\v_j)=\min_{i=1}^n(v(\v_i)$ and set $g_i=\v_i$ if $i\neq j$,
$g_j=\sigma$. We have $W_v=W_v(g_1,\ldots ,g_n)$. Also,
$v(\d^mg_i/\d t_v^m)\ge v(g_i)-m$ for $m\ge 0$, whence $$
v(W_v)\ge \sum_{i=1}^nv(g_i)-{n\choose
2}=v(\sigma)-v(\v_{j_v})+\sum_{i=1}^nv(\v_i)-{n\choose 2}. $$
Recalling that $v(\v_{j_v})=\min_{i=1}^nv(\v_i)$  we then obtain
$$ {n\choose 2}(2g-2)=\sum_vv(W_v)\ge \sum_{v\in
S}(v(\sigma)-\min_{i=1}^nv(\v_i))+\sum_{i=1}^n(\sum_{v\in
S}v(\v_i))-{n\choose 2}\# S. $$ Finally, for $i\le r$ all zeros
and poles of $\v_i$ are contained in $S$ so $\sum_{v\in
S}v(\v_i)=0$ for $i\le r$. For $i>r$ at least the poles of $\v_i$
are contained in $S$, so $\sum_{v\in S}v(\v_i)\ge -\deg(\v_i)$ for
$i> r$. Inserting this in the last displayed inequality yields the
sought result. \CVD\bigskip
 
The proof of Theorem 1 is based on two simple, though crucial,
points. The first one is embodied in the proof of the  following
statement, which is actually a weak form of Theorem 1.\medskip
 
\noindent{\bf Proposition 2}. {\it There exists a (computable)
function $B_1$ on $\N$ such that if $g,h\in\C[x]$ are non-constant
polynomials and if $g(h(x))$ has at most $l$ terms, then $h(x)$
may be written as a ratio of two polynomials  having  each at most
$B_1(l)$ terms.}
\medskip
 
\noindent{\it Proof}. It plainly suffices to construct the
function $B_1$ assuming that $h(x)$ is not of the shape $ax^n+b$.
Then, putting $d:=\deg g$ and $m:=\deg f$, $d$ is a divisor of
$m$,  $\deg h=m/d$ and by  Theorem A  we have $d\le 2l(l+1)$.
 
For $l=1$ we may take $B_1(1)=2$: in fact, if $g(h(x))$ is a
monomial $ax^m$, then $g$ cannot have two distinct roots and must
be therefore  of the shape $b(x-\xi)^n$. Then $h(x)=\xi+\eta
x^{m/n}$ ($b\eta^n=a$) has at most two terms.
 
We now argue by induction, supposing that $B_1$ has been suitably
defined on $\{1,\ldots ,l\}$.

We write $y=1/x$ and $f(x)=ax^{m}(1+b_1y^{n_1}+\ldots
+b_{l}y^{n_l})=ax^m\tilde f(y)$, say,  where $n_0:=0<n_1<\ldots
<n_l\le m$. We may suppose that $f$ has exactly $l+1$ terms, so
$ab_1\cdots b_l\neq 0$.
 
 From the equation $g(h)=f$ we may write  the Puiseux expansion for $h=h(x)$ at
$x=f=\infty$: $$
h(x)=c_{-1}f(x)^{1/d}+c_0+c_1f(x)^{-1/d}+\ldots,\qquad
c_{-1},c_0,c_1,\ldots \in\C,\eqno(1) $$ for a suitable choice of
the $d$-th root $f(x)^{1/d}$, where the $c_j$ depend only on $g$;
this identity is valid in $\C((y))$. We expand the various powers
of this $d$-th root
  as
$f(x)^{s/d}=a^{s/d}x^{ms/d}(1+b_1y^{n_1}+\ldots
+b_{l}y^{n_l})^{s/d}$, using the multinomial theorem for  the sum
on the right: $$ \tilde f(y)^{s/d}=(1+b_1y^{n_1}+\ldots
+b_{l}y^{n_l})^{s/d}=\sum_{h_1,\ldots
,h_l}c_{s,d,\h}b_1^{h_1}\cdots b_l^{h_l}y^{h_1n_1+\ldots
+h_ln_l},\eqno(2) $$ where $\h:=(h_1,\ldots ,h_l)$ runs through
$\N^l$ and where the $c_{s,d,\h}$ are certain universal
coefficients.

Factoring $g(x)$ we see that, since $x^{m-n_l}||f(x)=g(h(x))$,
there exists a  root $\xi$ of $g$, of multiplicity $d_0\le d$,
  such that $x^{m-n_l}||(h(x)-\xi)^{d_0}$.
Let us then write  $\tilde h(y)=x^{-m/d}(h(x)-\xi)\in\C[y]$. It
will suffice to prove the conclusion for  $\tilde h$ in place of
$h$.
 
We have $\deg\tilde h={m\over d}-{m-n_l\over d_0}\le {n_l\over d}$
(since $d_0\le d$ and $n_l\le m$). Also, subtracting $\xi$ from
both sides of (1) and dividing by $x^{m/d}$ we obtain, in the ring
$\C[[y]]$, for certain $\gamma_{-1},\gamma_0,\gamma_1,\ldots
\in\C$, $$ \tilde h(y)=\gamma_{-1}\tilde
f(y)^{1/d}+\gamma_0y^{m/d}+\gamma_1y^{2m/d}\tilde
f(y)^{-1/d}+\ldots,\eqno(3) $$
 
We note that  since $\tilde h(y)$ is a polynomial of degree $\le
n_l/d$, formula  (3) shows that it is the sum of  the terms on the
right of (2), with $s=1$, for which $h_1n_1+\ldots +h_ln_l\le
n_l/d$, plus possibly $\gamma_0y^{m/d}$. The number of such terms
is $\ll (n_l/n_1)^l$, hence, if we knew that $n_1>\epsilon_l n_l$
for some fixed $\epsilon_l>0$ we could easily establish the
conclusion of Theorem 1. This lower bound for $n_1$ isn't of
course guaranteed, but nevertheless we shall show that we can
somewhat reduce to this case.\medskip
 
We fix an integer $p$, $0\le p\le l-1$, and we write
$\delta_p(y)=1+b_1y^{n_1}+\ldots +b_py^{n_p}$ (so
$\delta_0(y)=1$).
 
Our  main task  will be now to establish that: {\it  if $n_{p+1}$
is not much smaller than $n_l$ (i.e. $\gg_l n_l$),
  then either we obtain the sought representation or  $n_p$ is as
   well not much smaller
than  $n_l$}. We shall then conclude by backward induction on
$p=l-1,l-2,\ldots $.\medskip

To take advantage of the fact that $n_{p+1}$ may be possibly
``large" we expand $\tilde f(y)^{s/d}$ in a slightly different
way, namely writing $$ \tilde
f(y)^{s/d}=\delta_p(y)^{s/d}\left(1+{b_{p+1}y^{n_{p+1}}+\ldots
+b_ly^{n_l}\over \delta_p(y)}\right)^{s/d} $$ and using the
multinomial theorem for the root of the sum on the right. In this
way each of the summands $\gamma_sy^{(1-s)m/d}\tilde f(y)^{s/d}$,
$s=-1,0,1,\ldots $,  on the right side of (3) will be expressed
(again in the ring $\C[[y]]$) as an infinite sum of terms of the
shape $$ c\delta_p(y)^{{s\over d}-k}y^{{(1-s)m\over
d}+h_1n_{p+1}+\ldots +h_{l-p}n_l},\qquad h_1+\ldots
+h_{l-p}=k,\quad c\in\C,\eqno(4) $$ for varying integers
$h_1,\ldots ,h_{l-p}\in\N$ and suitable constants $c=c(h_1,\ldots
,h_{l-p},s)$. Therefore $\tilde h(y)$ will be likewise expressed:
note in fact that $\delta_p(y)^{1/d}\in\C[[y]]^*$ and that the
infinite sum  converges in $\C[[y]]$  since $(1-s)m/d\to +\infty$
as $s\to -\infty$.
 
 We now consider all the terms of the shape (4)
such that the exponent of $y$ is $\le 2n_l$, i.e. ${(1-s)m\over
d}+h_1n_{p+1}+\ldots +h_{l-p}n_l\le 2n_l$. Clearly for this we
must have $s\ge 1-2d$ and $\max h_i\le 2n_l/n_{p+1}$. Hence the
number $L$ of such terms is bounded by a certain function  of $d$
and of $n_l/n_{p+1}$: for our present  purposes  we may take for
instance the rough estimate $L\le (2d+1)(1+(2n_l/n_{p+1}))^l$.
 
Denoting by $t_1,\ldots ,t_L$ such terms, we have in the ring
$\C[[y]]$, $$ \tilde h(y)=t_1+\ldots +t_L+O(y^{2n_l}).\eqno(5) $$
If $t_1,\ldots ,t_L$ are linearly dependent over $\C$, we may use
a linear relation to replace some $t_i$ by a linear combination of
the others. Hence, replacing $L$ with a possibly smaller number
  and changing if necessary the $t_i$ with suitable
constant multiples of themselves, we may assume that the $t_i$ in
(5) are linearly independent over $\C$ and that they are still of
the shape (4).\medskip

With the purpose of  applying Proposition 1,  we proceed to define
the relevant objects which appear in that statement.  We define
$K$ as the function field $\C(y,\delta_p(y)^{1/d})$. We readily
find $2g-2\le  dn_p$ for the genus. We let $n=L+1$, $\v_i:=-t_{i}$
for $i=1,\ldots ,L$, $\v_{L+1}:=\tilde h(y)$, so in fact $\v_i\in
K$ for all $i$. Also, $\sigma=\tilde h(y)-t_1-\ldots -t_L$. We
further let $r=L$ and we define $S$  as the set of zeros/poles of
$\v_1,\ldots ,\v_L$ together with the poles of $\v_{L+1}$. Now,
from (4) we see that $\v_1,\ldots ,\v_L$ have zeros/poles at most
at the places of $K$ above $0$ or $\infty$ of $\C(y)$ or above the
roots of $\delta_p(y)$, while $\v_{L+1}=\tilde h(y)$ has a pole
only at the places of $K$ above $\infty$.  This gives at most
$d(2+n_p)$ places in $S$.
\medskip
 
We now distinguish between two alternatives.
 
{\it First case}. This occurs when  $\v_1,\ldots ,\v_n$ are
linearly dependent over $\C$. In a relation of linear dependence
$\v_{L+1}$ must appear because we are assuming that $\v_1,\ldots
,\v_L$ are independent. Then we may express $\v_{L+1}=\tilde h(y)$
as a linear combination of at most $L$ terms of the shape (4). Let
$e:=[K:\C(y)]$, so $e$ is a certain divisor of $d$, in fact the
least integer such that $\delta_p(y)^{e}$ is a $d$-th power in
$\C(y)$.\note{In fact we must have $e=1$ in this case, but we
won't need this.} We may then write the said linear relation as
$$
\tilde h(y)=\sum_{j=0}^{e-1}\delta_p(y)^{j/d}\Lambda_j,
$$ where
$\delta_p(y)^{j/d}\Lambda_s$ is the sum of the terms of the shape
(4) in the linear relation, for which $s\equiv j\pmod {e}$; in
particular, $\Lambda_j\in\C(y)$. We deduce that $\tilde
h(y)=\Lambda_0$.
 
Note that $\delta_p(y)^{e/d}$ is a certain polynomial $\eta_p(y)$
such that $\eta_p(y)^{d/e}=\delta_p(y)$ has at most $p+1\le l$
terms. By the inductive assumption $\eta_p(y)$ can be written as a
ratio of two polynomials with at most $B_1(l)$ terms. Also,
$\Lambda_0$ is a sum of at most $L$ terms of the shape $$
c\eta_p(y)^{{s-kd\over e}}y^{{(1-s)m\over d}+h_1n_{p+1}+\ldots
+h_{l-p}n_l},\qquad e|(s,d),\quad |s|\le 2d-1, \quad \max (h_i)\le
2n_l/n_{p+1}. $$ In particular, since $k=\sum h_i\le
2ln_l/n_{p+1}$, since $d\le 2l(l+1)$ and since $L\le
(2d+1)(1+(n_l/n_{p+1}))^l$, $\tilde h(y)=\Lambda_0$   may be
written as a ratio of two polynomials each with $\le B_2(l,
n_l/n_{p+1})$ terms, where $B_2(l,u)$ is a function which may be
easily estimated in terms of $B_1(l)$ and of $u\ge 1$.

\medskip
{\it Second case}. Let us now analyze the remaining possibility,
i.e. that $\v_1,\ldots ,\v_n$ are linearly independent over $\C$.
In this case the conclusion of Proposition 1 holds.
 
The meaning of  (5) is that $v_0(\sigma)\ge 2n_lv_0(y)\ge 2n_l$
for some place $v_0$ of $K$ above the zero of $\C(y)$, so $v_0\in
S$. We clearly have $\min_{i=1}^nv_0(\v_i)\le v_0(\v_{L+1})=0$,
because $\tilde h(0)\neq 0$.

Since $v(\sigma)\ge \min_{i=1}^nv(\v_i)$  for all places $v$ of
$K$  and since $\deg_K\tilde h(y)\le d\deg\tilde h\le n_l$,
Proposition 1 yields $$ 2n_l\le {L+1\choose 2}(\# S+dn_p)+n_l\le
{(L+1)^2\over 2}(\# S+dn_p)+n_l. $$

 We have seen that $\# S\le d(2+n_p)$,  hence this
inequality becomes $$ n_l\le (L+1)^2d(1+n_p)\le
16^{l+1}d^3(n_l/n_{p+1})^{2l}(1+n_p),\eqno(6) $$ where we have
used our previous estimate for $L$. \medskip
 
Now, suppose that the first alternative never occurs, for
$p=l-1,l-2,\ldots ,0$. Then (6) is always true. For $p\ge 1$ it
gives $$ (n_l/n_p)\le 16^{l+2}d^3(n_l/n_{p+1})^{2l}. $$ Hence,
since $n_l/n_{p+1}=1$ for $p=l-1$, we obtain by backward induction
that $n_l/n_1$ is bounded in terms of $l$ only (recall $d\le
2l(l+1)$). We already noted that this suffices, but we may also
apply (6) for $p=0$ to get that $n_l$  is bounded only in terms of
$l$. Hence the degree of $\tilde h(y)$  and {\it a fortiori} the
number of its terms are  bounded by a (computable) function of
$l$, and we are done.
 
Therefore we may assume that the first alternative sometimes
occurs, and we denote by $q\ge 0$ the last such occurrence. Then
for $p>q$ the second alternative must hold, so we have (6) for
$p>q$. As before, inductively  we may  use this to show  that
$n_l/n_{q+1}$ is bounded by a function of $l$ only. Also,  since
the first alternative occurs for $p=q$, the previous argument
yields that $\tilde h(y)$ can be written as a ratio of two
polynomials whose number of terms is bounded by $B_2(l,
n_l/n_{q+1})$; but this is in turn bounded by a function only of
$l$, concluding finally the proof of Proposition 2. \CVD\bigskip
 
  To deduce
Theorem 1 from Proposition 2 we have only to show that $h(x)$ is
not just a ``ratio of polynomials with boundedly many terms", but
that itself has boundedly many terms; the examples $(x^n-1)/(x-1)$
show that this is not automatic. In our case, this will follow
from a description equivalent to Theorem 2, which we state
as:\medskip
 
\noindent{\bf Theorem $2^*$}. {\it Let $l$ be a positive integer
and write $f(x)=a_1x^{m_1}+\ldots +a_lx^{m_l}=g(h(x))\in \C[x]$,
where $\deg g,\deg h>1$ and  where $h(x)$ is not of the shape
$ax^n+b$.
 
Then $\deg g\le 2l(l+1)$ and $h(x)$ has at most $B=B(l)$ terms.
 
Further, there are finitely many algebraic varieties $\V_j\subset
\A^{l+2l(l+1)+1+B}$ (defined over $\Q$) and subgroups $\Lambda_j$
of $\Z^{l+B}$, $j=1,\ldots ,J=J(l)$,  such that for some $j\in
\{1,\ldots ,J\}$ the vector of coefficients of $f,g,h$ lies in
$\V_j$ and the vector of exponents of $x$ in $f,h$ lies in
$\Lambda_j$.
 
Conversely, if these vectors lie in $\V_j,\Lambda_j$ then the
equation $f=g\circ h$ holds.
 
Finally, one may effectively find $J$, defining equations for the
$\V_j$ and generators for the $\Lambda_j$.}\medskip

\noindent{\it Proof of Theorems 1,2$^*$}. Let $f(x)=g(h(x))$ where
$f(x)$ has at most $l$ terms and where $h(x)$ is not of the shape
$ax^n+b$. Then, from Theorem A and Proposition 2 it follows that
$\deg g\le 2l(l+1)=:\ell$ and that $h(x)$ is a ratio
$h_1(x)/h_2(x)$ where $h_1,h_2\in\C[x]$ have each at most
$B=B_1(l)$ terms. We may then write $$
f(x)=\sum_{i=1}^la_ix^{m_i},\quad g(x)=\sum_{j=0}^\ell
b_jx^j,\quad h_r(x)=\sum_{k=1}^Bc_{rk}x^{n_{rk}},\ r=1,2, $$ which
yields $$
(\sum_{k=1}^Bc_{2k}x^{n_{2k}})^\ell(\sum_{i=1}^la_ix^{m_i})-\sum_{j=0}^\ell
b_j(\sum_{k=1}^Bc_{1k}x^{n_{1k}})^j(\sum_{k=1}^Bc_{2k}x^{n_{2k}})^{\ell
-j}=0. $$ Expanding everything  we obtain the vanishing of a sum
of terms each of the shape $\gamma x^\mu$, where:
 
(i) The occurring degrees $\mu$ are certain explicitly given
linear combinations of $m_i,n_{1k},n_{2k}$ ($i=1,\ldots ,l$,
$k=1,\ldots ,B$) with coefficients in $\N$, bounded by $\ell$.
 
(ii) The occurring coefficients $\gamma$ are certain explicit
monomials (over $\Q$) in the $a_i,b_j,c_{1k},c_{2k}$ ($i=1,\ldots
,l$, $j=1,\ldots ,\ell$,  $k=1,\ldots ,B$), the set of these
monomials depending only on $l$.\medskip
 
We now group together all the terms $\gamma x^\mu$ having equal
degree $\mu$. This gives a partition of the terms, the possible
partitions being finite in number.
 
For each such partition, the various equalities between the
degrees gives (in view of (i)) a linear system with integral
coefficients, among the $m_i,n_{1k},n_{2k}$. Note that by (i) all
the  systems so obtained may be written down, and their number is
bounded only in terms of $l$. Forgetting the fact that the
$m_i,n_{rk}$ are non-negative, each system admits a parametric
solution of the shape $$ m_i=\sum_{j=1}^p\alpha_{ij}u_j,\quad
n_{rk}=\sum_{j=1}^p\beta_{rkj}u_j,\quad r=1,2,\eqno(7) $$ with
certain computable integers $p$, $\alpha_{ij},\beta_{rkj}$
depending only on the linear system and bounded only in terms of
$l$, where the $u_j$ may take any integer values; that is, for
arbitrary values of the $u_j$ the degrees are equal in groups
according to the partition (and conversely).
 
After grouping all the terms according to the partition, we equate
to zero all the corresponding coefficients. This gives an
algebraic system in the $a_i,b_j,c_{1k},c_{2k}$, defining a
certain affine algebraic variety (over $\Q$, possibly reducible).
These equations have only finitely many possibilities which can be
enumerated, and their  number is bounded only in terms of
$l$.\medskip
 
Plainly, each relevant equality $f(x)=g(h(x))$  produces then  a
point in one of these varieties and an integral solution of the
corresponding linear system. \note{Note also that {\it a priori}
one equation $f=g\circ h$ may give points and solutions in more
than one way; this is because certain equalities between degrees
do not exclude further equalities.}\medskip
 
Conversely, each point in a relevant variety together with an
integer solution of the corresponding linear system gives an
equality $f(x)=g(h(x))$, except for the fact that $f,h$ may now be
Laurent polynomials (that is, polynomials in $x,x^{-1}$). If we
restrict to non-negative integer solutions we obtain polynomial
equalities. \note{It may be easily seen that the non-negative
solutions of an integral linear system may be parametrized as
linear combinations with non-negative coefficients of a finite
system of generators.} \medskip
 
Now comes the second main point of the proof of Theorem 1. Take
any given equation $f(x)=g(h(x))$ of the shape in question. The
coefficients of $f,g,h_1,h_2$ will then give a certain point $P_0$
in one of the above varieties. And the degrees in $f,h_1,h_2$ will
give a solution $\w_0$  of the corresponding linear system,
expressed  by a parameter vector $\u_0$ as in (7). However, by
construction, if we keep the point $P_0$ fixed and vary the
solution $\w_0$ to any solution $\w$ of the same system, given by
the parameter $\u$, we shall   obtain another equation
$f^\u(x)=g(h_1^\u(x)/h_2^\u(x))$, with the same coefficients but
generally different degrees.
 
Let us exploit what this means. The degrees occurring in $f^\u$
are certain fixed linear combinations of $u_1,\ldots ,u_p$, with
integer coefficients depending on the  $\alpha_{ij}$ in (7), and
an analogous fact holds  for $h_1^\u,h_2^\u$. In other words, by
(7) we can write $$ f^\u(x)=F(x^{u_1},\ldots ,x^{u_p}),\qquad
h_r^\u(x)=H_r(x^{u_1},\ldots ,x^{u_p}),\qquad r=1,2, $$ for
certain Laurent polynomials $F,H_1,H_2\in \C[z_1^{\pm 1},\ldots
,z_p^{\pm 1}]$. Note that the degrees of the terms occurring in
$F,H_1,H_2$ depend only on the coefficients of the linear system
(which are bounded in terms of $l$) whereas the coefficients of
$F,H_1,H_2$ depend only on the linear system and  on the point
$P_0$, both of which are  fixed in the present discussion.

Since $f^\u=g(h_1^\u/h_2^\u)$ holds for all $\u\in\Z^p$, we deduce
that $F=g(H_1/H_2)$. But then $H_1/H_2$ is integral over the
integrally closed ring $\C[z_1^{\pm 1},\ldots ,z_p^{\pm 1}]$ and
therefore
  $H_1/H_2\in\C[z_1^{\pm 1},\ldots ,z_p^{\pm 1}]$. Note now that
  $\deg(H_1/H_2)\le \deg
H_1+\deg H_2$. Hence the number of terms of the Laurent polynomial
$H_1/H_2$ is bounded by $(1+2\deg H_1+2\deg H_2)^p$. But $p,\deg
H_1, \deg H_2$ depend only on the linear system, and are therefore
bounded only in terms of $l$.
 
This also shows that the number of terms of
$h^\u(x)=h_1^\u(x)/h_2^\u(x)$ is bounded by a function only of
$l$, and the same holds for $\u=\u_0$, proving Theorem 1.
 
Finally, Theorem 2$^*$ is obtained from Theorem 1 just by
repeating the opening arguments  of this proof; we  have only to
forget about $h_1(x),h_2(x)$ and keep only $h(x)$, which, as we
now know, has at most $B_3(l)$ terms.\CVD
\medskip
 
\noindent{\it Proof of Theorem 2}. Theorem 2 is a mere rephrasing
of Theorem 2$^*$. In fact, by construction, starting from the
variety $\V_j$ of Theorem 2$^*$ we may obtain $F_j,G_j,H_j$ as in
Theorem 2, just by setting $z_i:=x^{u_i}$ and by taking the
coefficients to be  the coordinate functions on $\V_j$.

\bigskip\medskip
 
\centerline{\title References}\medskip
 
\item{[BM]} - D. Brownawell and D. Masser, Vanishing sums in function fields, {\it
Math. Proc. Camb. Phil. Soc} {\bf 100} (1986), 427-434.\smallskip
 
\item{[E]} - P. Erd\"os, On the number of terms of the square of a polynomial, {\it
Niew Arch. Wiskunde} {\bf 23} (1949), 63-65.\smallskip
 
\item{[S]} - A. Schinzel, On the number of terms of a power of a polynomial, {\it Acta
Arith.}, {\bf XLIX} (1987), 55-70.\smallskip
 
\item{[S2]} - A. Schinzel, {\it Polynomials with special regard to reducibility},
Encyclopedia of Mathematics and its applications, vol. {\bf 77},
Cambridge Univ. Press, 2000.\smallskip
 
\item{[Z]} - U. Zannier, On the number of terms of a composite polynomial, preprint,
2006.\smallskip
 
\item{[Z2]} - U. Zannier, Some remarks on the $S$-unit equation in function fields,
{\it Acta Arith.} {\bf  LXIV} (1993), 87-98.\smallskip
 
\vfill
 
Umberto Zannier
 
Scuola Normale Superiore
 
Piazza dei Cavalieri 7
 
56126 Pisa - ITALY
 
email:  {\it u.zannier@sns.it}

\end